\documentclass{article}
\usepackage{graphicx}

\usepackage{amssymb}
\usepackage{bm}
\usepackage[toc,page]{appendix}

\begin{document}


\title{Isogeometric analysis for turbulent flow
\footnote{This work has been supported by Technology agency of the Czech Republic through the
project TA03011157 "Innovative techniques for improving utility qualities of water turbines
with the help of shape optimization based on modern methods of geometric modelling."}}

\maketitle

\begin{center}
\author{\large Bohum\'{i}r Bastl$^{a,b}$},
\author{\large Marek Brandner$^{a,b}$},
\author{\large Ji\v{r}\'{i} Egermaier$^{a,b}$},\\
\author{\large Krist\'{y}na Mich\'{a}lkov\'{a}$^{a,b}$},
\author{\large Eva Turnerov\'{a}$^{a,b}$}
\end{center}

\begin{center}
\emph{$^a$European Centre of Excellence New Technologies for the Information Society, University
of West Bohemia in Pilsen, Univerzitn\'{i} 22, 306 14 Plze\v{n}, Czech Republic}

\emph{$^b$University of West Bohemia, Department of Mathematics, Univerzitn\'{i} 22, 306 14 Plze\v{n},
Czech Republic}
\end{center}



\begin{abstract}
The article is devoted to the simulation of viscous incompressible turbulent fluid flow based on solving the Reynolds averaged Navier-Stokes (RANS) equations with different $k-\omega$ models. The isogeometrical approach is used for the discretization based on the Galerkin method. Primary goal of using isogeometric analysis is to be always geometrically exact, independent of the discretization, and to avoid a time-consuming generation of meshes of computational domains. For higher Reynolds numbers, we use stabilization SUPG technique in equations for $k$ and $\omega.$ The solutions are compared with the standard benchmark example of turbulent flow over a backward facing step.
\end{abstract}

\textbf{Keywords: } Isogeometric analysis; turbulent flow; $k-\omega$; Navier-Stokes equations;



\section{Introduction}

The main goal of our work is to propose and implement the numerical model for solving turbulent flow based on the Galerkin approach with NURBS basis functions. We show that it is possible to effectively solve turbulent flow described by the Navier-Stokes equations with $k-\omega$ turbulent model by the isogeometric manner.

The objectives of isogeometric analysis based on NURBS (non-uniform rational B-splines) are to generalize and improve finite element analysis. It means to
provide more accurate modelling of geometries and to exactly represent shapes
such as circles, cylinders, ellipsoids, etc. Due to exact geometry at the coarsest
level of discretization it is possible to eliminate geometrical errors. It also much simplifies mesh refinement of industrial geometries by eliminating communication with the CAD description of geometry. Further refinement of the mesh or increasing the order of basis functions are very simple, efficient and robust. At the same time, isogeometric analysis has many features in common with finite element analysis. For example the isoparametric concept in which dependent variables and the geometry share the same basis functions. Then the mesh, and
the corresponding basis, can be refined and order-elevated while maintaining the original exact geometry. Then the isogeometric methodology can be an useful tool for computational fluid dynamics, in particular, turbulent flows.

High Reynolds number turbulent flows are important in many applications. Turbulent 
flows involve multiscale space and time-developing flow physics. The dynamic of all relevant scales of the flow described by the Navier-Stokes equations can be solved by the direct numerical simulation (DNS) approach, which is too expensive for most practical flows. Therefore the most common approach is the Reynolds-Averaged Navier-Stokes (RANS), which simulates the mean flow and effects of all turbulent scales. The efficient intermediate approach is large eddy simulation (LES), which can simulate significant flow unsteadiness that RANS cannot handle. In the absence of universal turbulence theory there exist many developments and improvements of the schemes including the empirical information.

\section{Navier-Stokes equations}

The model of viscous flow of an incompressible Newtonian fluid can be described by the Navier-Stokes equations in the common form
\begin{equation}\label{NNS}
\begin{array}{rcl}
\displaystyle \frac{\partial \bm u}{\partial t} + \nabla p + \bm u \cdot \nabla \bm u - \nu \Delta \bm u & = & \bm f, \qquad \textrm{in } \Omega \times \langle 0,T\rangle,\\
\nabla \cdot \bm u & = & 0, \qquad \textrm{in } \Omega \times \langle 0,T\rangle,
\end{array}
\end{equation}
where $\Omega \subset \mathbb R^d$ (dimension $d = 1,2,3$) is the computational domain, $T > 0$ is the final time, $\bm u = \bm u(\bm x,t)$ is the vector function describing flow velocity, $p = p(\bm x,t)$ is the pressure function, $\nu$ describes kinematic viscosity and $f$ additional body forces acting on the fluid. We do not assume only very small Reynolds numbers, but there are still some "limits" for which this model gives reasonable solution or it is necessary to use very fine discretization. The initial-boundary value problem is considered as the system (\ref{NNS}) together with a suitable initial conditions and the following boundary conditions
\begin{equation}
\begin{array}{rclll}
\bm u & = & \bm w & \textrm{ on } \partial \Omega_D \times \langle 0,T\rangle& \textrm{(Dirichlet condition)},\\
\displaystyle \nu \frac{\partial \bm u}{\partial \bm n} - \bm n p & = & \bm 0 & \textrm{ on } \partial \Omega_N \times \langle 0,T\rangle & \textrm{(Neumann condition)}.
\end{array}
\end{equation}
If the velocity is specified everywhere on the boundary, then the pressure solution is only unique up to a hydrostatic constant. 

The Navier-Stokes equations describe turbulent incompressible flow. This flow contains many eddies of different sizes which are changing in time. The numerical methods for solving turbulent models are divided into the following categories:

\begin{itemize}
\item \textbf{Direct numerical simulation (DNS)}

The average flux and all turbulent fluctuations are computed. It means, that we use FEM (or FVM etc.) method directly to solve the Navier-Stokes equations. It is necessary to use a very fine mesh to compute small turbulent fluctuations of the flow and a short time step for non-stationary problem. Therefore this approach is very computationally expensive. 

\item \textbf{Reynolds Averaged Navier Stokes (RANS)}

This approach simulates only average flux and effects of this flux to the flow. It uses the time averaged Navier-Stokes equations. The special term appears in the equations which is approximated by the appropriate approaches. The most common approaches are $k-\varepsilon$ or $k-\omega$ models. RANS is very often used in practice. 

\item \textbf{Large Eddy Simulation (LES)}

LES simulates behaviour of large eddies which is realized by averaging of the Navier-Stokes equations in space dimension thus the small eddies are not considered. The behaviour of small eddies is described by the so called subgrid scale model, which can be for example computed by RANS.
\end{itemize}

\section{RANS models}

RANS is the most common model for solving the Navier-Stokes equations including turbulence. It is based on decomposition of the solution into the time-averaged value and fluctuation value. In two dimensions, the solution $\bm u = \bm u(x_1,x_2,t)$ and $p = p(x_1,x_2,t)$ is decomposed by the following way 
\begin{equation}
\bm u\left(x_1, x_2, t\right) = \bar{\bm u}\left(x_1, x_2\right) + \bm u'\left(x_1, x_2, t\right),
\label{eq:dekomp1}
\end{equation}
\begin{equation}
p\left(x_1, x_2, t\right) = \bar{p}\left(x_1, x_2\right) + p'\left(x_1, x_2, t\right),
\label{eq:dekomp2}
\end{equation}
where $\bar{\bm u}, \bar{p}$ are time-averaged values and $\bm u', p'$ are 
fluctuation ones. Substituting (\ref{eq:dekomp1}) and (\ref{eq:dekomp2}) to the system (\ref{NNS}) we arrive at (see \cite{Versteeg} for details)
\begin{equation}
\frac{\partial \bar{\bm{u}}}{\partial t} + \bar{\bm{u}}\cdot\nabla \bar{\bm{u}} = -\nabla \bar{p} + \nu\Delta \bar{\bm{u}} - \overline{\bm{u}'\cdot\nabla \bm{u}'},
\label{eq:rans1}
\end{equation}
\begin{equation}
\nabla\cdot \bar{\bm{u}} = 0,
\label{eq:rans2}
\end{equation}
The solution are the functions $\bar{\bm u}, \bar{p}$, the fluctuation values are not determined. The equations contain the unknown term
\begin{equation}
\overline{\bm{u}'\cdot\nabla \bm{u}'}.
\end{equation}
This term is approximated by the following relation (see \cite{Davidson2} for details)
\begin{equation}
- \overline{\bm{u}'\cdot\nabla \bm{u}'} = \nabla\cdot\left(\nu_{T}\left(\nabla\bar{\bm{u}} + (\nabla\bar{\bm{u}})^{T}\right) - \frac{2}{3}k\bm I\right),
\label{eq:stress_aprox_matrix}
\end{equation}
where $\bm I$ is identity matrix and (in 2D case)
$$
\nabla\bar{\bm{u}} + (\nabla\bar{\bm{u}})^{T} =
\left[\begin{array}{cc}
	\frac{\partial\bar{u}_{1}}{\partial x_1} & \frac{\partial\bar{u}_{1}}{\partial x_2}  \\
	\frac{\partial\bar{u}_{2}}{\partial x_1} & \frac{\partial\bar{u}_{2}}{\partial x_2}  \\
\end{array}\right] + \left[\begin{array}{cc}
	\frac{\partial\bar{u}_{1}}{\partial x_1} & \frac{\partial\bar{u}_{2}}{\partial x_1}  \\
	\frac{\partial\bar{u}_{1}}{\partial x_2} & \frac{\partial\bar{u}_{2}}{\partial x_2}  \\
\end{array}\right] = \left[\begin{array}{cc}
	\frac{\partial\bar{u}_{1}}{\partial x_1} + \frac{\partial\bar{u}_{1}}{\partial x_1} & \frac{\partial\bar{u}_{1}}{\partial x_2} + \frac{\partial\bar{u}_{2}}{\partial x_1}  \\
	\frac{\partial\bar{u}_{1}}{\partial x_2} + \frac{\partial\bar{u}_{2}}{\partial x_1} & \frac{\partial\bar{u}_{2}}{\partial x_2} + \frac{\partial\bar{u}_{2}}{\partial x_2}  \\
\end{array}\right].
$$
Function $\nu_T$ is the so-called turbulent viscosity defined by
\begin{equation}\label{eq13}
\nu_{T} = \frac{k}{\omega},
\end{equation}
where $k$ is the turbulent kinetic energy and $\omega$ is the specific dissipation. Our goal is to determine these quantities $k$ and $\omega$ in order to approximate (\ref{eq:stress_aprox_matrix}). Several approaches can be applied but the common way is $k-\omega$ model which is described in the following part. 

Substituting (\ref{eq:stress_aprox_matrix}) to (\ref{eq:rans1}) RANS equations can be written in the form
\begin{equation}
\frac{\partial \bar{\bm{u}}}{\partial t} + \bar{\bm{u}}\cdot\nabla \bar{\bm{u}} = -\nabla \bar{p} + \underbrace{\nabla\cdot\left[\left(\nu + \nu_{T}\right)\nabla\bar{\bm{u}}\right]}_{\nu\Delta\bar{\bm{u}} + \nabla\cdot(\nu_{T}\nabla \bar{\bm{u}})} + \nabla\cdot(\nu_{T}\left(\nabla\bar{\bm{u}}\right)^{T}) - \frac{2}{3}\nabla k,
\label{eq:rans_1}
\end{equation}
\begin{equation}
\nabla\cdot \bar{\bm{u}} = 0.
\label{eq:rans_2}
\end{equation}

\subsection{Basic $k-\omega$ model}\label{sec31}

This model adds two extra equations to the RANS system, i.e. the transport equation for the turbulent kinetic energy $k$ and the equation for the specific dissipation $\omega.$

Transport equations for kinetic energy and specific dissipation have the form
\begin{equation}\label{eq14}
\frac{\partial k}{\partial t} + \bar{\bm{u}}\cdot\nabla k = \nabla\left[\left(\frac{\nu_T}{\sigma_k} + \nu\right)\nabla k\right] + \nu_Tf - C_{\mu}k\omega,
\end{equation}
\begin{equation}\label{eq15}
\frac{\partial\omega}{\partial t} + \bar{\bm{u}}\cdot\nabla\omega = \nabla[(\sigma_{\omega}\nu_T + \nu)\nabla\omega] + C_{\omega1}\frac{\omega}{k}\nu_Tf - C_{\omega2}\omega^2,
\end{equation}
where $\bar{\bm{u}} = (\bar{u}_1,\bar{u}_2), f = \left(\frac{\partial\bar{u}_1}{\partial x_2} + \frac{\partial\bar{u}_2}{\partial x_1}\right)^2 + 2\left(\frac{\partial\bar{u}_1}{\partial x_1}\right)^2 + 2\left(\frac{\partial\bar{u}_2}{\partial x_2}\right)^2 = \frac{1}{2}|\nabla\bar{\bm{u}} + \nabla\bar{\bm{u}}^T|^2,$ and $\sigma_{\omega}, \sigma_k, C_{\mu},C_{\omega1}$ and $C_{\omega2}$ are constants with the standard choice of values: $\sigma_{\omega} = 0.5, \sigma_k = 2.0, C_{\mu} = 0.09, C_{\omega1} = 0.52$ and $C_{\omega2} = 0.072.$ For details about $k-\omega$ model see \cite{Celik,Davidson1,Petit}.

The solution of system (\ref{eq14}) and (\ref{eq15}) defines the turbulent viscosity $\nu_T$ by the relation (\ref{eq13}).

For simplicity the equations (\ref{eq14}) and (\ref{eq15}) can be separated with the help of (\ref{eq13})
\begin{equation}\label{eq16}
\frac{\partial k}{\partial t} + \bar{\bm{u}}\cdot\nabla k = \nabla\cdot\left[\left(\frac{\nu_T}{\sigma_k} + \nu\right)\nabla k\right] + \nu_Tf - C_{\mu}\frac{k^2}{\nu_T},
\end{equation}
\begin{equation}\label{eq17}
\frac{\partial\omega}{\partial t} + \bar{\bm{u}}\cdot\nabla\omega = \nabla\cdot[(\sigma_{\omega}\nu_T + \nu)\nabla\omega] + C_{\omega1}f - C_{\omega2}\omega^2.
\end{equation}

\subsection{Low Reynolds number model (LRN)}\label{sec32}

One of the common problems in turbulent modelling is computing turbulent flow influenced by the adjacent wall. There is a boundary layer, where the velocity changes from the no-slip condition at the wall to its free stream value. The standard method of solving this problem is to apply a very fine mesh close to the wall. This is the so called integration method, which necessitates an LRN type of turbulence model. On the other hand, this method does not need any wall-function approximation which is a common tool for obtaining near wall values in other methods. The region under wall influence diminishes in the case of higher
Reynolds numbers.

We shortly describe two turbulent LRN models, which we use for numerical experiments.

\subsubsection{Wilcox 1993}

This model is in detail described in \cite{Bredberg}. The equations for $k$ and $\omega$ are defined as
\begin{equation}\label{eq18}
\frac{\partial k}{\partial t} = P_k - \beta^*k\omega + \nabla\cdot\left[(\nu + \sigma^*\nu_t)\nabla k\right],
\end{equation}
\begin{equation}\label{eq19}
\frac{\partial\omega}{\partial t} = \alpha\frac{\omega}{k}P_k - \beta\omega^2 + \nabla\cdot\left[(\nu + \sigma^*\nu_t)\nabla\omega\right],
\end{equation}
where
$$
\nu_t = \alpha^*\frac{k}{\omega}, \qquad \alpha^* = \frac{\alpha^*_0 + \frac{Re_T}{R_k}}{1 + \frac{Re_T}{R_k}}, \qquad \beta^* = \frac{9}{100}\frac{\frac{5}{18} + \left(\frac{Re_T}{R_{\beta}}\right)^4}{1 + \left(\frac{Re_T}{R_{\beta}}\right)^4},
$$
$$
\alpha = \frac{5}{9\alpha^*}\frac{\alpha_0 + \frac{Re_T}{R_{\omega}}}{1 + \frac{Re_T}{R_{\omega}}}, \quad Re_T = \frac{k}{\nu\omega}, \quad P_k = 2\nu_tS_{ij}\frac{\partial u_i}{\partial x_j}, \quad S_{ij} = \frac{1}{2}\left(\frac{\partial u_i}{\partial x_j} + \frac{\partial u_j}{\partial x_i}\right),
$$
$$
\alpha_0 = \frac{1}{10}, \alpha_0^* = \frac{\beta}{3}, \beta = \frac{3}{40}, \sigma^* = \frac{1}{2}, R_{\beta} = 8, R_k = \frac{27}{10}, R_{\omega} = 6,
$$
where $Re_T$ is turbulent Reynolds number. The turbulent viscosity is defined as $\nu_T = \nu_t$ in this model.

\subsubsection{Low Reynolds number version of Wilcox (2006) $k-\omega$ two-equation model}

This model is in detail described in \cite{Wilcox}. The equations for $k$ and $\omega$ are defined as
\begin{equation}\label{eq20}
\frac{\partial k}{\partial t} + \bar{\bm{u}}\cdot\nabla k = P - \beta^*k\omega + \nabla\cdot\!\left[(\nu + \sigma_k\nu_t)\nabla k\right],
\end{equation}
\begin{equation}\label{eq21}
\frac{\partial\omega}{\partial t} + \bar{\bm{u}}\cdot\!\nabla\omega = \frac{\gamma\omega}{k}P - \beta\omega^2 + \nabla\cdot\left[(\nu + \sigma_{\omega}\nu_t)\nabla\omega\right] + \frac{\sigma_{d}}{\omega}\nabla k\cdot\nabla\omega,
\end{equation}
where
$$
P = \tau_{ij}\frac{\partial u_i}{\partial x_j}, \qquad \tau_{ij} = \nu_t\left(2S_{ij} - \frac{2}{3}\frac{\partial u_k}{\partial x_k}\delta_{ij}\right) - \frac{2}{3}k\delta_{ij},
$$
$$
\beta^* = 0.09\left(\frac{100\frac{\beta_0}{27} + \left(\frac{Re_T}{R_{\beta}}\right)^4}{1 + \left(\frac{Re_T}{R_{\beta}}\right)^4}\right), \qquad \gamma = \frac{13}{25}\left(\frac{\alpha_0 + \frac{Re_T}{R_{\omega}}}{1 + \frac{Re_T}{R_{\omega}}}\right)(\alpha^*)^{-1},
$$
$$
\sigma_d = \left\{\begin{array}{lr}
0, & \textrm{ for } \frac{\partial k}{\partial x_j}\frac{\partial \omega}{\partial x_j} \leq 0,\\
\frac{1}{8}, & \textrm{ for } \frac{\partial k}{\partial x_j}\frac{\partial \omega}{\partial x_j} > 0,\end{array}\right.
\quad Re_T = \frac{k}{\nu\omega}, \quad \alpha^* = \frac{\alpha_0^* + \frac{Re_T}{R_k}}{1 + \frac{Re_T}{R_k}}, \quad \alpha_0^* = \frac{\beta_0}{3},
$$
$$
R_{\beta} = 8, \quad R_k = 6, \quad R_{\omega} = 2.61, \quad \alpha_0 = \frac{1}{9}, \quad \beta_0 = 0.0708.
$$
As mentioned above, solving the two equation model leads to the approximation of the Reynolds stresses (\ref{eq:stress_aprox_matrix}). However, using LRN Wilcox $k-\omega$ model (2006), the turbulent viscosity $\nu_{T}$ is not defined by the expression (\ref{eq13}) again, but the modified expression is used instead
$$
\nu_{T} = \alpha^{*}\frac{k}{\hat{\omega}},
$$
where
$$
\hat{\omega} = \max\left[\omega,C_{\lim}\sqrt{\frac{2S_{ij}S_{ij}}{\frac{\beta^*}{\alpha^*}}}\right], \qquad C_{\lim} = \frac{7}{8}.
$$

\section{Numerical model}

In this section, we describe some numerical aspects and techniques on which the model
is based. 

\subsection{NURBS approximation}

\begin{figure}[ht!]
\begin{center}
\begin{tabular}{cc}
$T=(0,1,2,\ldots,7), p=1$ & $T=(0,1,2,\ldots,7), p=3$\\
\includegraphics[width=4cm, height=2cm]{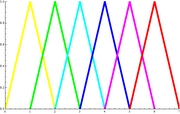} & \includegraphics[width=4cm, height=2cm]{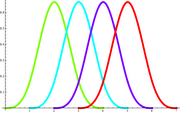}\\
$T = (0, 0, 0, 0, 1, 2, 3, 4, 4, 4, 4)$ & $T = (0, 0, 0, 1, 2, 2, 3, 3, 3)$\\
\includegraphics[width=4.5cm, height=2cm]{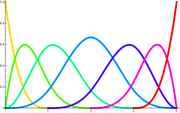} & \includegraphics[width=4.5cm, height=2cm]{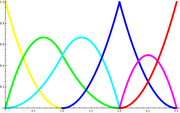}
\end{tabular}
\end{center}
\caption{B-spline basis functions.}\label{obrBBF}
\end{figure}

NURBS surface of degree $p$, $q$ is determined by a control net $\bm{P}$ (of control points $P_{i,j}$, $i=0,\ldots,n$, $j=0,\ldots,m$), weights $w_{i,j}$ of these control points and two knot vectors $U=(u_0,\ldots,u_{n+p+1})$, $V=(v_0,\ldots,v_{m+q+1})$ and is given by a parametrization
\begin{equation}
S(u,v) = \frac{\sum_{i=0}^n \sum_{j=0}^m w_{i,j} P_{i,j} N_{i,p}(u) M_{j,q}(v)}{\sum_{i=0}^n \sum_{j=0}^m w_{i,j} N_{i,p}(u) M_{j,q}(v)} = \sum_{i=0}^n \sum_{j=0}^m P_{i,j} R_{i,j}(u,v).
\end{equation}
B-spline basis functions $N_{i,p}(u)$ and $M_{j,q}(v)$ are determined by knot vectors $U$ and $V$ and degrees $p$ and $q$, respectively, by a formula (for $N_{i,p}(u)$, $M_{j,q}(v)$ is constructed by the similar way)
\begin{eqnarray}\label{bspliebazedef}
    N_{i,0}(u) & = & \left\{
                    \begin{array}{ll}
                    1 & u_i \leq t < u_{i+1}, \\
                    0 & \mbox{otherwise},
                    \end{array}
                    \right. \nonumber \\
    N_{i,p}(u) & = & \frac{u-u_i}{u_{i+p}-u_i} N_{i,p-1}(u) + \frac{u_{i+p+1}-u}{u_{i+p+1}-u_{i+1}} N_{i+1,p-1}(u). 
    \end{eqnarray}
Knot vector is a non-decreasing sequence of real numbers which determines the distribution of a parameter on the corresponding curve/surface. B-spline basis functions (see Figure \ref{obrBBF}) of degree $p$ are $C^{p-1}$-continuous in general. Knot repeated $k$ times in the knot vector decreases the continuity of B-spline basis functions by $k-1.$ Support of B-spline basis functions is local -- it is nonzero only on the interval $[t_i,t_{i+p+1}]$ in the parameter space and each B-spline basis function is non-negative, i.e., $N_{i,p}(t)\geq 0, \forall t.$

\subsection{Galerkin approach}

Let $V$ be a velocity solution space and $V_0$ be the corresponding space of test functions, i.e.,
\begin{equation}
\begin{array}{rcl}
V & = & \{ \bm u \in H^1(\Omega)^d | \bm u = \bm w \textrm{ on } \partial \Omega_D\},\\
V_0 & = & \{ \bm v \in H^1(\Omega)^d | \bm v = \bm 0 \textrm{ on } \partial \Omega_D\}.\\
\end{array}
\end{equation}
Then a weak formulation of the stationary boundary value problem is: find $\bm u \in V$ and $p \in L_2(\Omega)$ such that
\begin{eqnarray}
\int_\Omega\bm v\cdot\frac{\partial\bm u}{\partial t} + \nu \int_\Omega \nabla \bm u : \nabla \bm v + \int_\Omega (\bm u \cdot \nabla \bm u) \bm v - \int_\Omega p \nabla \cdot \bm v & = & \int_\Omega \bm f \cdot \bm v, \forall \bm v \in V_0, \nonumber \\
\int_\Omega q \nabla \cdot \bm u & = & 0, \qquad\forall q \in L_2(\Omega).\nonumber
\end{eqnarray}
Further, we use Galerkin approach with general basis functions. The discrete weak solution is defined by finite dimensional spaces $V^h \subset V$, $W^h \subset L_2(\Omega)$ and their basis functions. Then we look for $\bm u_h \in V^h$ and $p_h \in W^h$ such that for all test functions $\bm v_h \in V_0^h \subset V_0$ and $q_h \in W^h$
\begin{eqnarray}
\int_\Omega\bm v_h\cdot\frac{\partial\bm u_h}{\partial t} + \nu\int\limits_{\Omega}\nabla\bm u_h:\nabla\bm v_h + \int_\Omega(\bm u_h\cdot\nabla\bm u_h)\bm v_h - \int\limits_{\Omega}p_h\nabla\cdot\bm v_h & = & \int\limits_{\Omega}\bm f\cdot\bm v_h,\nonumber\\
\int\limits_{\Omega}q_h\nabla\cdot\bm u_h & = & 0.
\end{eqnarray}
Isogeometric approach consists in taking the solution $\bm u_h$ as a linear combination of basis functions $R_i^u \in V^h$ and the solution $p_h$ is written as a linear combination of basis functions $R^p_i \in W^h,$ where $R_i^u$ and $R_i^p$ are NURBS description of a computational domain. In 2D, the solution has the form
\begin{equation}\label{LCBF}
\bm u_h = \sum_{i=1}^{n^u_{d}} (u_{1i}, u_{2i})^T R^u_i + \sum_{i=n^u_{d}+1}^{n^u_{v}} (u_{1i}^*, u_{2i}^*)^T R^u_i, \quad p_h = \sum_{i=1}^{n^p} p_i R_i^p,
\end{equation}
where $n^u_d$ is the number of points where the Dirichlet boundary condition is not defined. The discrete problem can be written (using implicit Euler method) in matrix form
\begin{eqnarray}
\left[
\begin{array}{ccc}
\bm A + \bm C & \bm 0 & -\bm B_1^T\\
\bm 0 & \bm A + \bm C & -\bm B_2^T\\
\bm B_1 & \bm B_2 & \bm 0
\end{array}
\right]
\left[
\begin{array}{c}
\bm u_1^{n+1}\\ \bm u_2^{n+1} \\ \bm p^{n+1}
\end{array}
\right] = \\
=
\left[\begin{array}{cc}
\bm C & \bm 0\\
\bm 0 & \bm C\\
\bm 0 & \bm 0
\end{array}\right]\left[\begin{array}{c}
\bm u_1^n\\
\bm u_2^n\end{array}\right]
 - \left[
\begin{array}{cc}
\bm A^* + \bm C^* & \bm 0 \\
\bm 0 & \bm A^* + \bm C^*\\
\bm B_1^* & \bm B_2^*
\end{array}
\right] \left[
\begin{array}{c}
\bm u_1^*\\
\bm u_2^*
\end{array}
\right],\nonumber
\end{eqnarray}
where
\begin{equation}
\begin{array}{lcllcl}
\bm A & = & \bigl[A_{ij}\bigr]_{1\leq i\leq n_d^u,1\leq j\leq n_d^u}, & \bm A^* & = &  \bigl[A_{ij}\bigr]_{1\leq i\leq n_d^u,n_d^u+1\leq j\leq n_v^u}, \\
\bm B_m & = &  \bigl[B_{mij}\bigr]_{1\leq i\leq n^p,1\leq j\leq n_d^u}, & \bm B_m^* & = &  \bigl[B_{mij}\bigr]_{1\leq i\leq n^p,n_d^u+1\leq j\leq n_v^u},\\
\bm C & = &  \bigl[C_{ij}\bigr]_{1\leq i\leq n_v^u,1\leq j\leq n_d^u}, & \bm C^* & = & \bigl[C_{ij}\bigr]_{1\leq i\leq n_d^u,n_d^u+1\leq j\leq n_v^u},
\end{array}
\end{equation}
\begin{equation}
\begin{array}{rcl}
A_{ij} & = & \displaystyle \nu\Delta t \int\limits_\Omega (\nabla R_i^u \cdot J^{-1}) \cdot (\nabla R_j^u \cdot J^{-1}) |\det J|,\\
B_{mij} & = & \displaystyle \Delta t\int\limits_\Omega R_i^p \bigl[ (\nabla R_j^u \cdot J^{-1}) \cdot \bm e_m\bigr] |\det J|,\\
C_{ij} & = & \displaystyle \int\limits_{\Omega} R_i^u R_j^u |\det J|,
\end{array}
\end{equation}
$m = \{1,2\}$ and $J$ is Jacobi matrix of a mapping from parametric domain to the computational domain. If $n^p > n^u_d$ then the system matrix has not the full rank - so we can choose the basis $R^u_i$ and $R^p_i$ arbitrary.

\subsection{Nonlinear iteration}

Because of non-linearity of Navier-Stokes equations it is necessary to solve the problem iteratively with linear problem in every step. One of the possibilities is to use the Picard's method. Here we mention the stationary problem for simplicity, the time derivative has no effect to the explanation of this method. 

First the non-linear residuum from weak formulation is computed by the values $\bm u^k$ and $p^k$ (for example by the solution of the Stokes problem). The residuum $R^k$ and $r^k$ satisfy
\begin{eqnarray}\label{Rk}
R^k &=& \int\limits_\Omega \bm f \cdot \bm v - \nu \int\limits_\Omega \nabla \bm u^k : \nabla \bm v - \int\limits_\Omega (\bm u^k \cdot \nabla \bm u^k) \bm v + \int\limits_\Omega p^k \nabla \cdot \bm v,\\ \label{rk}
r^k &=& - \int\limits_\Omega q \nabla \cdot \bm u^k.
\end{eqnarray}
We consider the exact solution $\bm u$ and $p$ as the sum of the solution of the current iteration and fluctuation of the solution $\bm u = \bm u^k + \delta\bm u^k$ and $p = p^k + \delta p^k.$ Substituting this solution to the weak formulation and using the equalities (\ref{Rk}) and (\ref{rk}) we have
\begin{eqnarray}\label{linuh}
R^k &=& \int\limits_\Omega \bm f \cdot \bm v - \nu \int\limits_\Omega \nabla \delta\bm u^k : \nabla \bm v - \int\limits_\Omega (\bm u^k \cdot \nabla \delta\bm u^k) \bm - \nonumber\\
& & - \int\limits_\Omega (\delta\bm u^k \cdot \nabla \bm u^k) \bm - \int\limits_\Omega (\delta\bm u^k \cdot \nabla \delta\bm u^k) \bm v + \int\limits_\Omega p^k \nabla \cdot \bm v,\\
r^k &=& - \int\limits_\Omega q \nabla \cdot \delta\bm u^k.\nonumber
\end{eqnarray}
We assume that the quadratic term $\int\limits_\Omega (\delta\bm u^k \cdot \nabla \delta\bm u^k) \bm v$ and also the linear term $\int\limits_\Omega (\delta\bm u^k \cdot \nabla \bm u^k) \bm v$ are sufficiently small and
we neglect them. We obtain the linear problem (\ref{linuh}) for fluctuations $\delta \bm u^k \in V_0$ and $\delta p^k \in L_2(\Omega).$ These define the following
step $\bm u^{k+1} = \bm u^k + \delta\bm u^k.$ Therefore we search for $\bm u^{k+1} \in V$ and $p^{k+1} \in L_2(\Omega)$ so that for all functions $\bm v \in V_0$ and $q \in L_2(\Omega)$ the following relation is valid
\begin{eqnarray}
\nu \int\limits_\Omega \nabla \bm u^{k+1} : \nabla \bm v + \int\limits_\Omega (\bm u^k \cdot \nabla \bm u^{k+1}) \bm v - \int\limits_\Omega p^{k+1} \nabla \cdot \bm v & = & \int\limits_\Omega \bm f \cdot \bm v,\nonumber\\
\int\limits_\Omega q \nabla \cdot \bm u^{k+1} & = & 0.
\end{eqnarray}
The solution of the Stokes problem is used as the initial condition for the iterative cycle.

We use Galerkin method and define the finite-dimensional spaces $V^h \subset V,$ $V_0^h \subset V_0, W^h \subset L_2(\Omega)$ and their basis functions. Find $\bm u_h \in V^h$ and $p_h \in W^h$ so that all functions $\bm v_h \in V_0^h$ a $q_h \in W^h$ satisfy
\begin{eqnarray}
\label{rce1}\nu \int\limits_\Omega \nabla \bm u_h^{k+1} : \nabla \bm v_h + \int\limits_\Omega (\bm u_h^k \cdot \nabla \bm u_h^{k+1}) \bm v_h - \int\limits_\Omega p_h^{k+1} \nabla \cdot \bm v_h & = & \int\limits_\Omega \bm f \cdot \bm v_h, \\
\label{rce2}\int\limits_\Omega q_h \nabla \cdot \bm u_h^{k+1} & = & 0.
\end{eqnarray}
The solution $\bm u_h^k$ and $p_h^k$ is written as linear combination of the basis functions (see (\ref{LCBF})) and it is substituted to (\ref{rce1}) and (\ref{rce2}). The sequence of the solutions $(\bm u_h^k, p_h^k) \in V^h \times W^h$ converges to the weak solution. The system has the matrix form:
\begin{eqnarray}
\left[
\begin{array}{ccc}
\bm A+\bm N(\bm u^k) & \bm 0 & -\bm B_1\\
\bm 0 & \bm A+\bm N(\bm u^k) & -\bm B_2\\
\bm B_1 & \bm B_2 & \bm 0
\end{array}
\right]
\left[
\begin{array}{c}
\bm u_1^{k+1}\\ \bm u_2^{k+1} \\ \bm p^{k+1}
\end{array}
\right] = \nonumber \\
= \left[
\begin{array}{c}
\bm f_1 - (\bm A^*+\bm N^*(\bm u^k)) \cdot \bm u_1^*\\
\bm f_2 - (\bm A^*+\bm N^*(\bm u^k)) \cdot \bm u_2^*\\
\bm B_1^* \cdot \bm u_1^* + \bm B_2^* \cdot \bm u_2^*
\end{array}
\right],
\end{eqnarray}
where
\begin{eqnarray*}
\bm A = \bigl[A_{ij}\bigr]_{1\leq i\leq n_d^u,1\leq j\leq n_d^u}, &\bm A^* =  \bigl[A_{ij}\bigr]_{1\leq i\leq n_d^u,n_d^u+1\leq j\leq n_v^u}, \\
\bm N(\bm u) = \bigl[N_{ij}(\bm u)\bigr]_{1\leq i\leq n_d^u,1\leq j\leq n_d^u}, &\bm N^*(\bm u) = \bigl[N_{ij}(\bm u)\bigr]_{1\leq i\leq n_d^u,n_d^u+1\leq j\leq n_v^u}, \\
\bm B_m = \bigl[B_{mij}\bigr]_{1\leq i\leq n^p,1\leq j\leq n_d^u},& \bm B_m^*  = \bigl[B_{mij}\bigr]_{1\leq i\leq n^p,n_d^u+1\leq j\leq n_v^u},
\end{eqnarray*}
\begin{equation}
\begin{array}{rcl}
A_{ij} & = & \displaystyle \nu \int\limits_\Omega (\nabla R_i^u \cdot J^{-1}) \cdot (\nabla R_j^u \cdot J^{-1}) |\det J|,\\
 N_{ij}(\bm u) &  = & \displaystyle  \int\limits_\Omega R_i^u \left[ \left( \sum_{l=1}^{n_v^u} (u_{1l}, u_{2l}) R_l^u\right) \cdot (\nabla R_j^u \cdot J^{-1})\right] |\det J|,\\
B_{mij} & = & \displaystyle \int\limits_\Omega R_i^p \bigl[ (\nabla R_j^u \cdot J^{-1}) \cdot \bm e_m\bigr] |\det J|.
\end{array}
\end{equation}

\subsection{SUPG - Streamline Upwind/Petrov-Galerkin}

Solving of the advection-diffusion equations can lead to numerical nonstability (for example the Navier-Stokes equations for high Reynolds numbers). One of the methods to reduce nonphysical oscillations is based on the construction of test function $\overline{\bm v}$ in special form (see for example \cite{John})
\begin{equation}
\overline{\bm  v} = \bm v + \tau_S \bm u \cdot \nabla \bm v,
\end{equation}
where
\begin{equation}
\tau_S = \frac{h}{2 \deg(\bm u)\|\bm u\|} \left(\coth P - \frac{1}{P}\right)
\end{equation}
and $h$ is the element diameter in the direction of $\bm u$ and $P = \frac{\|\bm u\|h}{2\nu}$ is the local P\'{e}clet number which determines whether the problem is locally convection dominated or diffusion dominated. In our test examples we use this SUPG stabilization method only for solving the $k-\omega$ equations, it is not applied to the solution of RANS equations.  

\section{Computational scheme (algorithm)}

This section is devoted to a more precise description of systems of linear algebraic equations arising from solving RANS equations with $k-\omega$ model via an isogeometric approach. Currently, these linear systems are solved with direct solver and we investigate suitable preconditioning strategies for preconditioned GMRES to be used. The implicit Euler method is used for the time discretization.

\subsection{Navier-Stokes equations}

We formulate weak solution for RANS system (\ref{eq:rans_1}), (\ref{eq:rans_2}). Let us consider function spaces $V \!=\! \left\{\bm{\varphi} \in H^{1}\left(\Omega\right): \bm{\varphi}|_{\partial\Omega} = \bm{w}\right\}$, $V_{0} \!=\! \left\{\bm{\psi} \in H^{1}\left(\Omega\right): \bm{\psi}|_{\partial\Omega} = \bm{0}\right\}.$ Then we search for $\bar{\bm u} \in V$, $\bar{p} \in L^{2}(\Omega)$ and $k \in L^{2}(\Omega)$ so that all test functions $\bm{v} \in V_{0}$ and $q \in L^{2}(\Omega)$ satisfy
$$
\int_{\Omega}\bar{\bm u}^{n+1}\cdot\bm v\mathrm{d}\Omega + \Delta t\nu\int_{\Omega}\nabla\bar{\bm u}^{n+1}\cdot\nabla\bm v\mathrm{d}\Omega - \Delta t\int_{\Omega}\overline{p}^{n+1}\cdot\nabla\bm v\mathrm{d}\Omega = 
$$
\begin{equation}\label{eq38}
= \int_{\Omega}\bar{\bm u}^n\cdot\bm v\mathrm{d}\Omega - \Delta t\int_{\Omega}(\bar{\bm u}^{n+1}\cdot\nabla\bar{\bm u}^{n+1})\bm v\mathrm{d}\Omega - 
\end{equation}
$$
- \Delta t\int_{\Omega}\nu_T^n\nabla\bar{\bm u}^{n+1}\cdot\nabla\bm v\mathrm{d}\Omega - \Delta t\int_{\Omega}\nu_T^n(\nabla\bar{\bm u}^{n+1})^T\cdot\nabla\bm v\mathrm{d}\Omega + \frac{2}{3}\Delta t\int_{\Omega}k^n\cdot\nabla\bm v\mathrm{d}\Omega,
$$
\begin{equation}\label{eq39}
\int_{\Omega}q\nabla\cdot\bar{\bm u}^{n+1}\mathrm{d}\Omega = 0.
\end{equation}
Now we consider finite dimensional spaces $V^{h} \subset V$, $V^{h}_{0} \subset V_{0}$ and $W^{h} \subset L^{2}(\Omega)$ with basis functions $R_{i}^{u}\in V^{h}$, $R_{i}^{p}$ and $R_{i}^{k}\in W^{h}$. The solution $\bar{\bm{u}}_{h} \in V^{h}$, $\bar{p}_{h} \in W^{h}$ and $k_{h} \in W^{h}$ have the forms of linear combinations of basis functions
\begin{equation}
\bar{\bm{u}}_{h}^{n+1} = \sum^{N}_{i = 1}(\bar{u}_{1i}^{n+1},\bar{u}_{2i}^{n+1})R_{i}^{u}, \quad \bar{\bm{u}}_{h}^{n} = \sum^{N}_{i = 1}(\bar{u}_{1i}^{n},\bar{u}_{2i}^{n})R_{i}^{u},
\end{equation}
\begin{equation}
\bar{p}_{h}^{n+1} = \sum^{N}_{i = 1}\bar{p}^{n+1}_{i} R_{i}^{p}, \quad \bar{p}_{h}^{n} = \sum^{N}_{i = 1}\bar{p}^{n}_{i} R_{i}^{p},
\end{equation}
\begin{equation}
k_{h}^{n} = \sum^{N}_{i = 1}k^{n}_{i} R_{i}^{k}.
\end{equation}
We consider $\bm v_h \in V_0^h.$ Substituting to (\ref{eq38}) and (\ref{eq39}) we obtain
\begin{eqnarray*}
\sum\limits_{i=1}^N(\bar{u}_{1i}^{n+1},\bar{u}_{2i}^{n+1})[\underbrace{\int_{\Omega}R_i^uR_j^u\mathrm{d}\Omega}_{C} + \Delta t\nu\underbrace{\int_{\Omega}\nabla R_i^u\nabla R_j^u\mathrm{d}\Omega}_{A} + \nonumber\\
+ \Delta t\underbrace{\int_{\Omega}\!\!\left(\sum\limits_{l=1}^N(\bar{u}_{1l}^{n},\bar{u}_{2l}^{n})R_l^u\right)\!\!\nabla R_i^uR_j^u\mathrm{d}\Omega}_{N(\bar{\bm u})} + \Delta t\underbrace{\int_{\Omega}\!\!\nu_T^n\nabla R_i^u\nabla R_j^u\mathrm{d}\Omega}_{E_1} + \Delta t(Q_1)_{i,j}] - \nonumber\\
- \sum\limits_{i=1}^N\bar{p}_i^{n+1}\Delta t\underbrace{\int_{\Omega}R_i^p\nabla R_j^u\mathrm{d}\Omega}_{B_1,B_2} = \sum\limits_{i=1}^N(\bar{u}_{1i}^{n},\bar{u}_{2i}^{n})[\underbrace{\int_{\Omega}R_i^uR_j^u\mathrm{d}\Omega}_{C}] - \Delta t(Q_2)_{i,j} + 
\end{eqnarray*}
\begin{equation}\label{eq45}
+ \frac{2}{3}\Delta t\sum\limits_{i=1}^Nk_i^n\underbrace{\int_{\Omega}R_i^k\nabla R_j^u\mathrm{d}\Omega}_{F_1,F_2},
\end{equation}
\begin{equation}\label{eq46}
\sum\limits_{i=1}^N(\bar{u}_{1i}^{n+1},\bar{u}_{2i}^{n+1})\underbrace{\int_{\Omega}R_j^p\nabla\cdot R_i^u\mathrm{d}\Omega}_{B_1,B_2} = 0,
\end{equation}
where
$$
(Q_{1})_{i,j} = \left[\begin{array}{c}
	\overbrace{\int_{\Omega}\nu_{T}^{n}\frac{\partial R_{i}^{u}}{\partial x_1}\frac{\partial R_{j}^{u}}{\partial x_1}\mathrm{d}\Omega}^{E_{2}} \\
	\underbrace{\int_{\Omega}\nu_{T}^{n}\frac{\partial R_{i}^{u}}{\partial x_2}\frac{\partial R_{j}^{u}}{\partial x_2}\mathrm{d}\Omega}_{E_{3}} 
\end{array}\right],
(Q_{2})_{i,j} = \left[\begin{array}{c}
	\sum^{N}_{i = 1}\bar{u}_{2i}^{n}\overbrace{\int_{\Omega}\nu_{T}^{n}\frac{\partial R_{i}^{u}}{\partial x_1}\frac{\partial R_{j}^{u}}{\partial x_2}\mathrm{d}\Omega}^{E_{4}} \\
	\sum^{N}_{i = 1}\bar{u}_{1i}^{n}\underbrace{\int_{\Omega}\nu_{T}^{n}\frac{\partial R_{i}^{u}}{\partial x_2}\frac{\partial R_{j}^{u}}{\partial x_1} \mathrm{d}\Omega}_{E_{5}}
\end{array}\right].
$$
The equations (\ref{eq45}), (\ref{eq46}) can be written in the matrix form
$$
\left[\begin{array}{ccc}
	\bm A_1 & \bm 0 & -\Delta t\bm B_{1}^{T} \\
	\bm 0 & \bm A_2 & -\Delta t\bm B_{2}^{T} \\
	\bm B_{1} & \bm B_{2} & \bm 0  
\end{array}\right]
\left[\begin{array}{c}
\bar{\bm{u}}_{1}^{n+1} \\
\bar{\bm{u}}_{2}^{n+1} \\
\bar{\bm{p}}^{n+1} 
\end{array}\right] =  
\frac{2\Delta t}{3}\left[\begin{array}{c}
\bm F_{1}^{T}\\
\bm F_{2}^{T}\\
0
\end{array}\right]\bm{k}^{n} + 
$$
\begin{equation}
+ \left[\begin{array}{ccc}
	\bm C & - \Delta t\bm E_{4} & \bm 0  \\
	- \Delta t\bm E_{5} & \bm C & \bm 0  \\
	\bm 0 & \bm 0 & \bm 0  
\end{array}\right]
\left[\begin{array}{c}
\bar{\bm{u}}_{1}^{n} \\
\bar{\bm{u}}_{2}^{n} \\
\bar{\bm{p}}^{n} 
\end{array}\right],
\label{ar:NS1}
\end{equation}
where
$$
\bm A_1 = \bm C + \Delta t\nu \bm A + \Delta t\bm N(\bar{\bm u}) + \Delta t\bm E_1 + \Delta t\bm E_2,
$$
$$
\bm A_2 = \bm C + \Delta t\nu \bm A + \Delta t\bm N(\bar{\bm u}) + \Delta t\bm E_1 + \Delta t\bm E_3
$$ 
and 
$$
\bar{\bm{u}}_{1}^{n+1} = \left(\bar{u}_{1 1}^{n+1}, \bar{u}_{1 2}^{n+1}, \ldots, \bar{u}_{1 N}^{n+1}\right), \bar{\bm{u}}_{1}^{n} = \left(\bar{u}_{1 1}^{n}, \bar{u}_{1 2}^{n}, \ldots, \bar{u}_{1 N}^{n}\right),
$$
$$
\bar{\bm{u}}_{2}^{n+1} = \left(\bar{u}_{2 1}^{n+1}, \bar{u}_{2 2}^{n+1}, \ldots, \bar{u}_{2 N}^{n+1}\right), \bar{\bm{u}}_{2}^{n} = \left(\bar{u}_{2 1}^{n}, \bar{u}_{2 2}^{n}, \ldots, \bar{u}_{2 N}^{n}\right),
$$
$$
\bar{\bm{p}}^{n+1} = \left(\bar{p}_{1}^{n+1}, \bar{p}_{2}^{n+1}, \ldots, \bar{p}_{N}^{n+1}\right), \bar{\bm{p}} ^{n} = \left(\bar{p}_{1}^{n}, \bar{p}_{2}^{n}, \ldots, \bar{p}_{N}^{n}\right),
$$
$$
\bm{k}^{n} = \left(k_{1}^{n}, k_{2}^{n}, \ldots, k_{N}^{n}\right).
$$

\subsection{$k-\omega$ turbulent model}

We consider $k-\omega$ model based on the equations (\ref{eq16}), (\ref{eq17}). Similar approach can be used for LRN models described in the section \ref{sec32}. If $k^n, \omega^n$ are approximations of unknown functions at the time step $t_n,$ then the approximations $k^{n+1}, \omega^{n+1}$ can be determined by the following relations (implicit Euler method is used, $\bar{\bm u}$ is supposed from previous time layer $t_n$)
\begin{equation}\label{eq48}
\frac{k^{n+1} - k^n}{\Delta t} + \bar{\bm u}^n\cdot\nabla k^{n+1} = \nabla\left[\left(\frac{\nu_T^n}{\sigma_k} + \nu\right)\nabla k^{n+1}\right] + \nu_T^nf^n - C_{\mu}\frac{(k^{n+1})^2}{\nu_T^n},
\end{equation}
\begin{equation}\label{eq49}
\frac{\omega^{n+1} - \omega^n}{\Delta t} + \bar{\bm u}^n\cdot\nabla \omega^{n+1} = \nabla[(\sigma_{\omega}\nu_T^n + \nu)\nabla \omega^{n+1}] + C_{\omega 1}f^n - C_{\omega 2}(\omega^{n+1})^2.
\end{equation}

Using the similar approach as in the previous section (see Appendix for details and definitions of matrices) we can formulate the following system
$$
\left[\begin{array}{cc}
\bm X_{11} & \bm 0\\
\bm 0 & \bm X_{22}\end{array}\right]\left[\begin{array}{c}
\bm k^{n+1}\\
\bm \omega^{n+1}\end{array}\right] =
$$
\begin{equation}
= \left[\begin{array}{cc}
\bm Y_{11} & \bm 0\\
\bm 0 & \bm Y_{22}\end{array}\right]\left[\begin{array}{c}
\bm k^n\\
\bm \omega^n \end{array}\right] + \left[\begin{array}{c}
\Delta t\bm F^k(\nu_T^n,\bm k^{n+1}) + \Delta t\nabla \bm F^k(\nu_T^n,\bm f^n)\\
\Delta tC_{\omega 1}\bm F^\omega(\bm f^n) + \Delta tC_{\omega 1}\nabla \bm F^\omega(\bm f^n)\end{array}\right],
\end{equation}
where
\begin{eqnarray*}
\bm X_{11} &=& \bm A^k + \Delta t\bm C^k(\nu_T^n) + \Delta tC_{\mu}\bm B_1^k(\nu_T^n,\bm k^{n+1}) + \nabla \bm A^k + \Delta t\nabla \bm C^k(\nu_T^n) + \\
& & + \Delta t\bm D^k + \Delta t\nabla \bm D^k + \Delta tC_\mu\nabla \bm B_1^k(\nu_T^n,\bm k^{n+1}),\\
\bm X_{22} &=& \bm A^\omega + \Delta t\bm C^\omega(\nu_T^n) + \Delta tC_{\omega 2}\bm B_1^\omega(\omega^{n+1}) + \nabla \bm A^\omega + \Delta t\nabla \bm C^\omega(\nu_T^n) + \\
& & + \Delta t\bm D^\omega + \Delta t\nabla \bm D^\omega + \Delta tC_{\omega 2}\nabla \bm B_1^\omega(\omega^{n+1}),\\
\bm Y_{11} &=& \bm A^k + \nabla\bm A^k,\\
\bm Y_{22} &=& \bm A^\omega + \nabla \bm A^\omega ,
\end{eqnarray*}
and $\bm k^{n+1} = (k_1^{n+1}, k_2^{n+1},\ldots,k_N^{n+1})^T, \bm k^n = (k_1^n, k_2^n,\ldots,k_N^n)^T$\\ and $\bm \omega^{n+1} = (\omega_1^{n+1}, \omega_2^{n+1},\ldots,\omega_N^{n+1})^T, \bm \omega^n = (\omega_1^n, \omega_2^n,\ldots,\omega_N^n)^T.$ 

\section{Numerical experiments}

\begin{figure}[ht!]
\begin{center}
\begin{tabular}{cc}
\includegraphics[width=10.7cm,height=2.7cm]{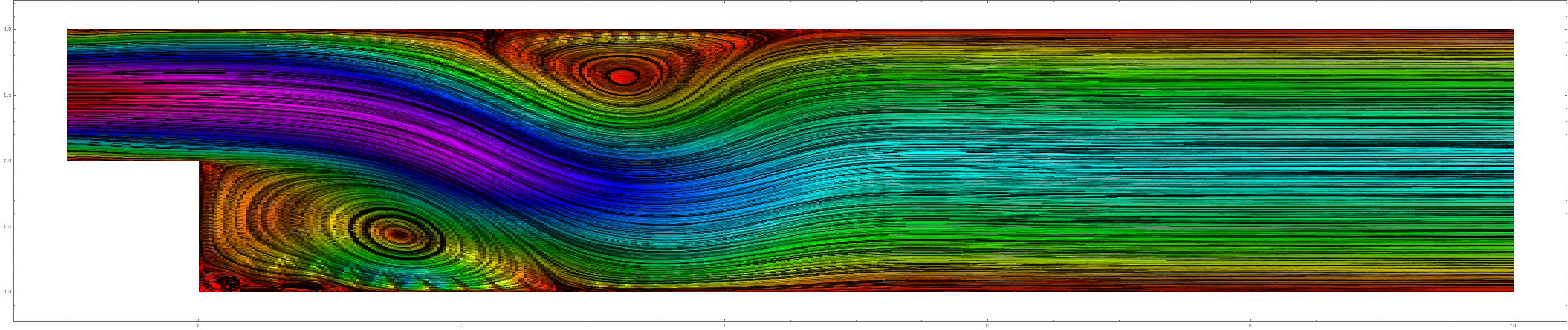} & \\
\includegraphics[width=10.7cm,height=2.7cm]{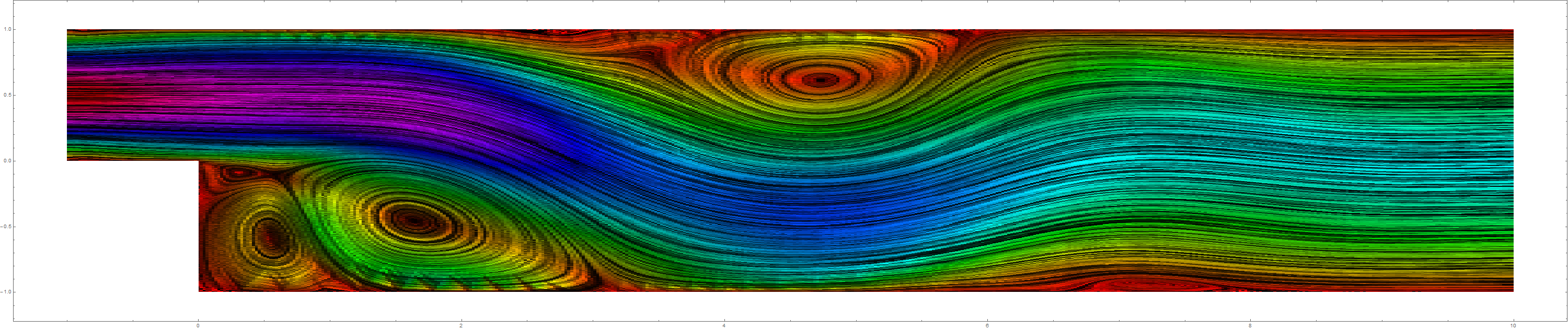} & \\
\includegraphics[width=10.7cm,height=2.7cm]{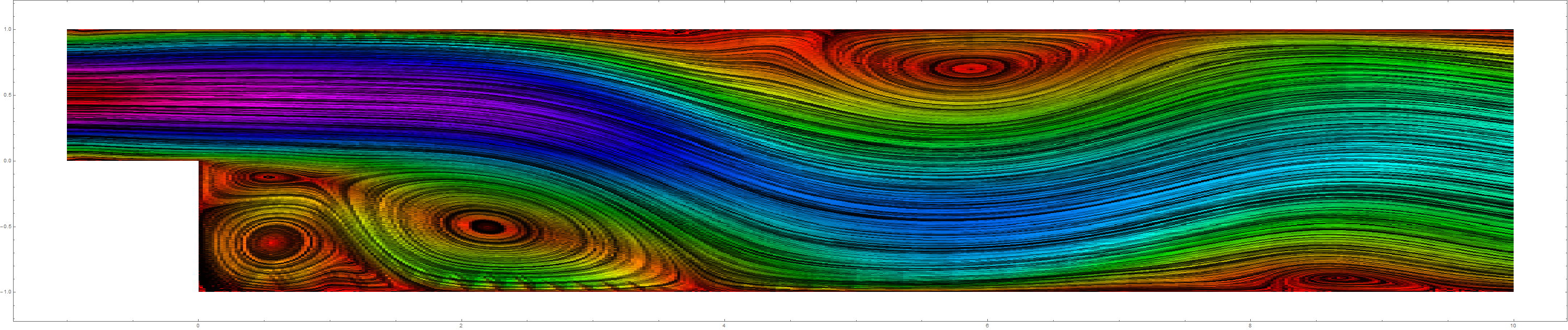} & \\
\includegraphics[width=10.7cm,height=2.7cm]{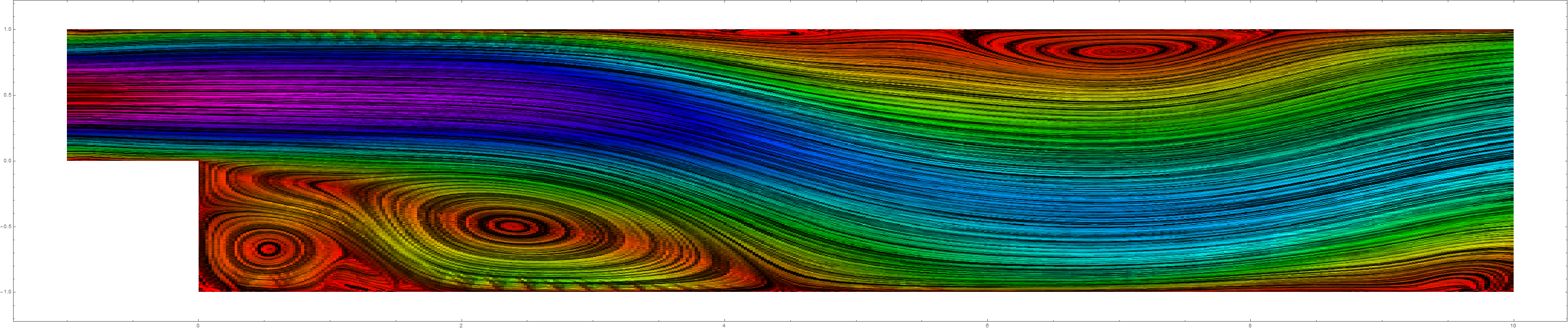} & \\
\includegraphics[width=10.7cm,height=2.7cm]{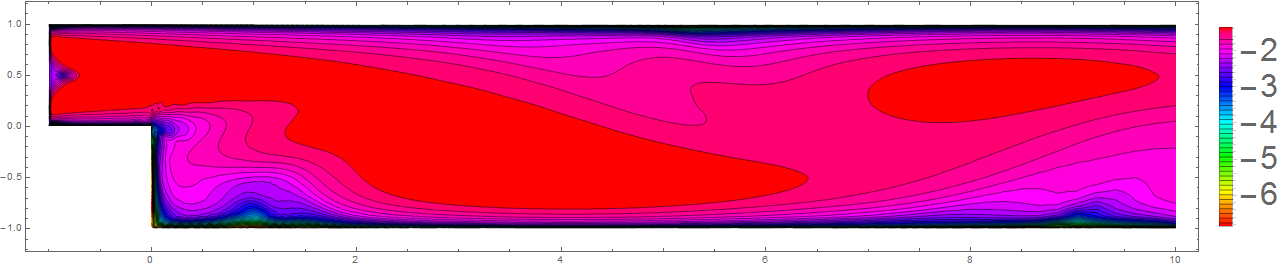}
\end{tabular}
\end{center}
\caption{The solution of turbulent flow based on turbulent model (\ref{eq16}), (\ref{eq17}). The upper pictures illustrate streamlines at the times $2.5s, 5s, 7.5s$ and $10s.$ The lower picture illustrates distribution of the turbulent viscosity (in logaritmic scale) at the time $10s.$}\label{obr:HRN}
\end{figure}

\begin{figure}[ht!]
\begin{center}
\begin{tabular}{cc}
\includegraphics[width=10.7cm,height=2.7cm]{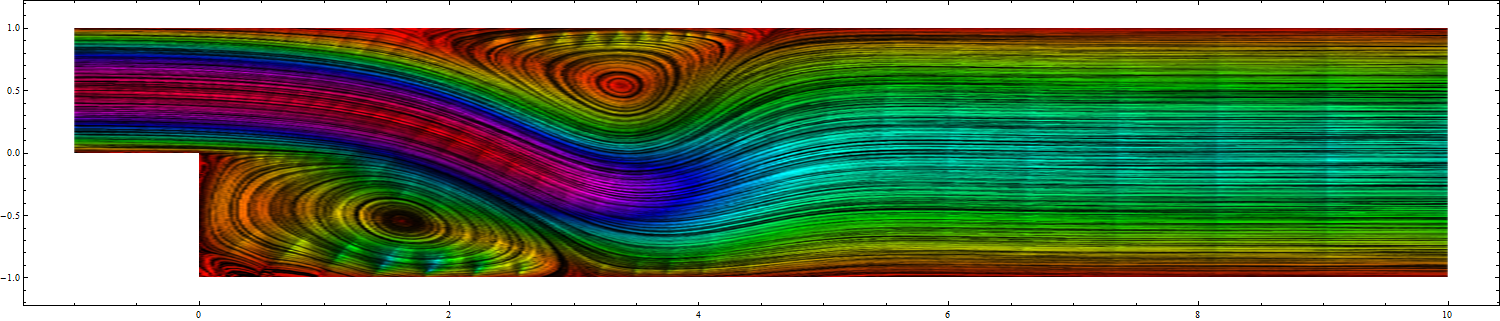} & \\
\includegraphics[width=10.7cm,height=2.7cm]{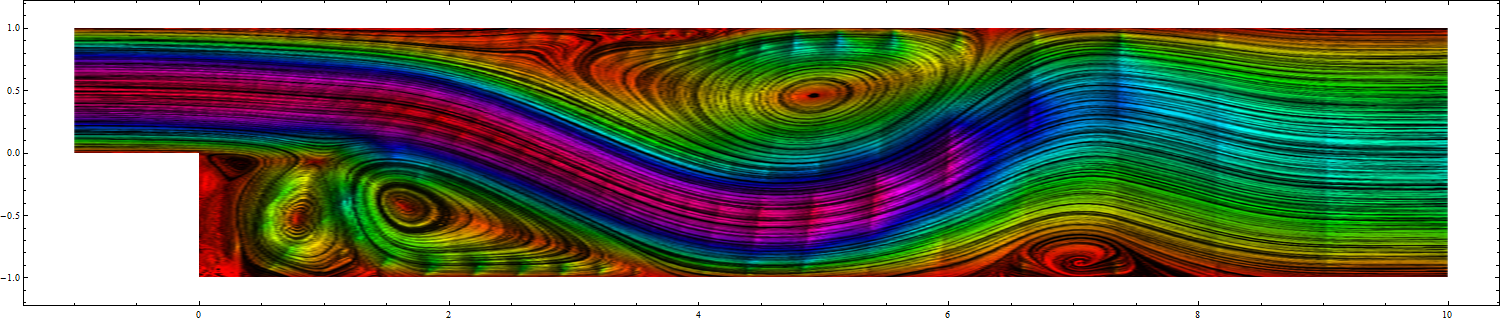} & \\
\includegraphics[width=10.7cm,height=2.7cm]{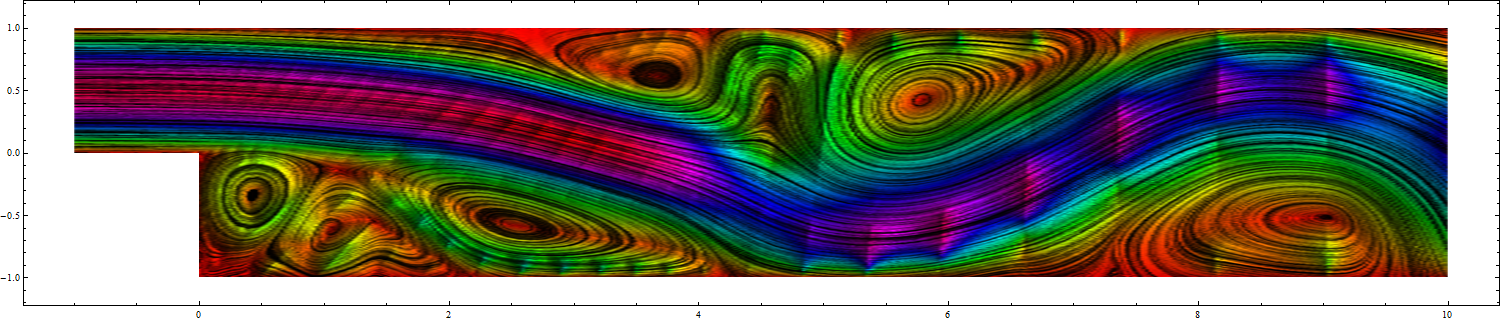} & \\
\includegraphics[width=10.7cm,height=2.7cm]{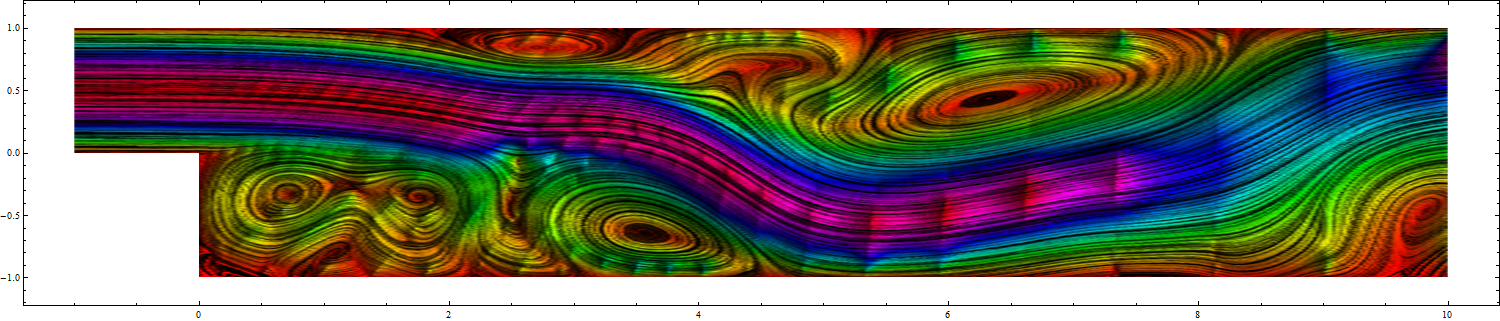} & \\
\includegraphics[width=10.7cm,height=2.7cm]{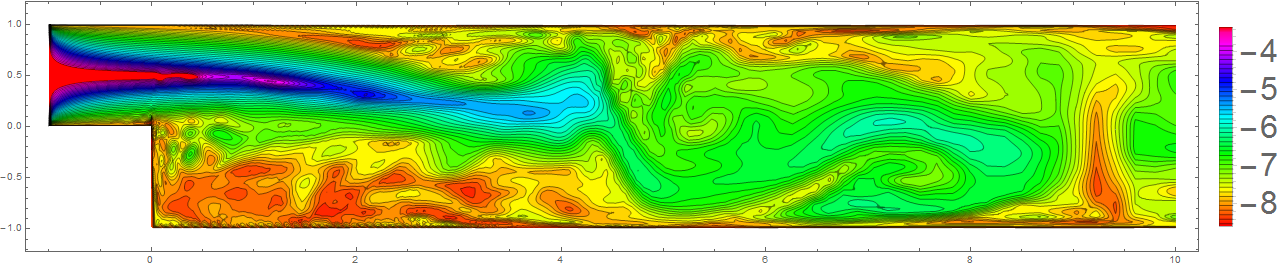}
\end{tabular}
\end{center}
\caption{The solution of turbulent flow based on turbulent model (\ref{eq18}), (\ref{eq19}). The upper pictures illustrate streamlines at the times $2.5s, 5s, 7.5s$ and $10s.$ The lower picture illustrates distribution of the turbulent viscosity (in logaritmic scale) at the time $10s.$}\label{obr:LRN1993}
\end{figure}

\begin{figure}[ht!]
\begin{center}
\begin{tabular}{cc}
\includegraphics[width=10.7cm,height=2.7cm]{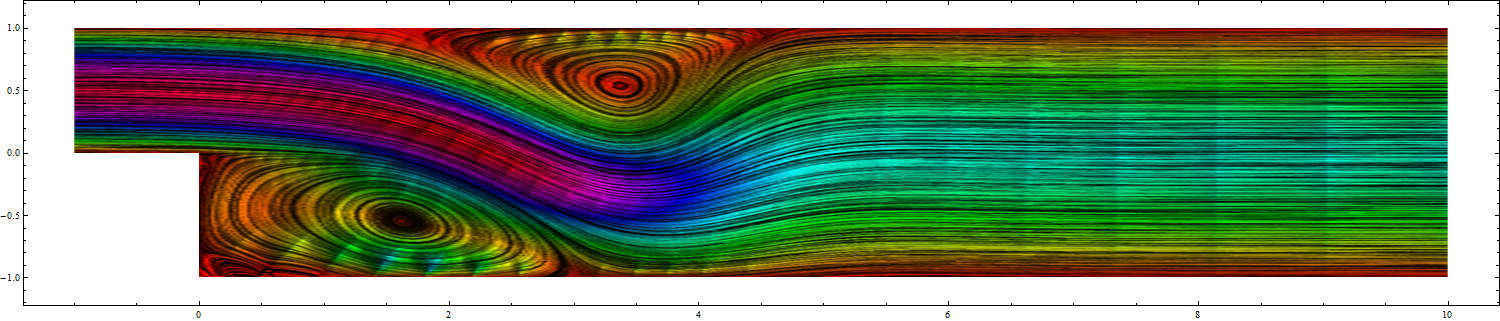} & \\
\includegraphics[width=10.7cm,height=2.7cm]{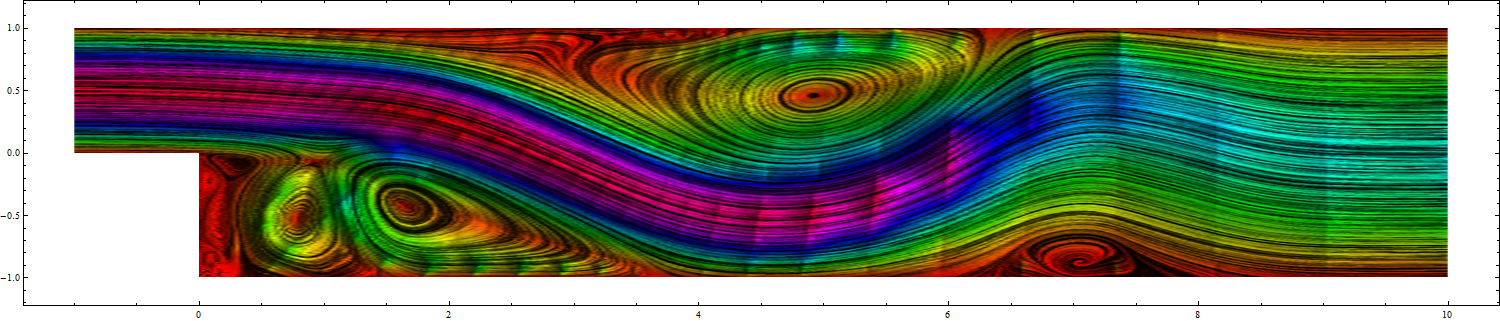} & \\
\includegraphics[width=10.7cm,height=2.7cm]{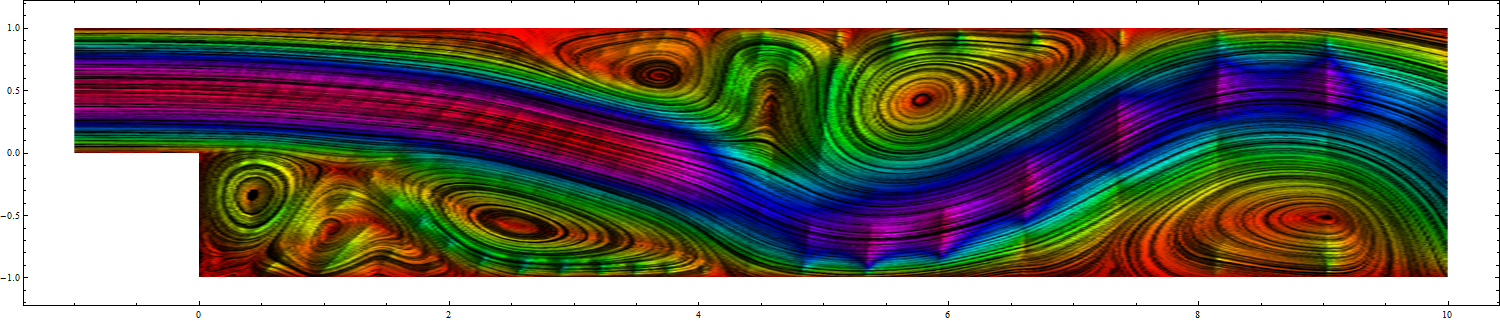} & \\
\includegraphics[width=10.7cm,height=2.7cm]{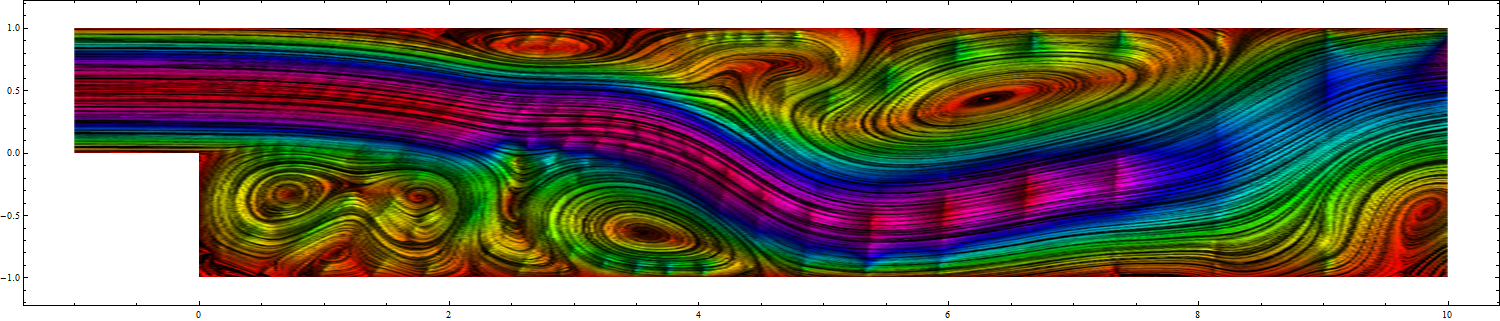} & \\
\includegraphics[width=10.7cm,height=2.7cm]{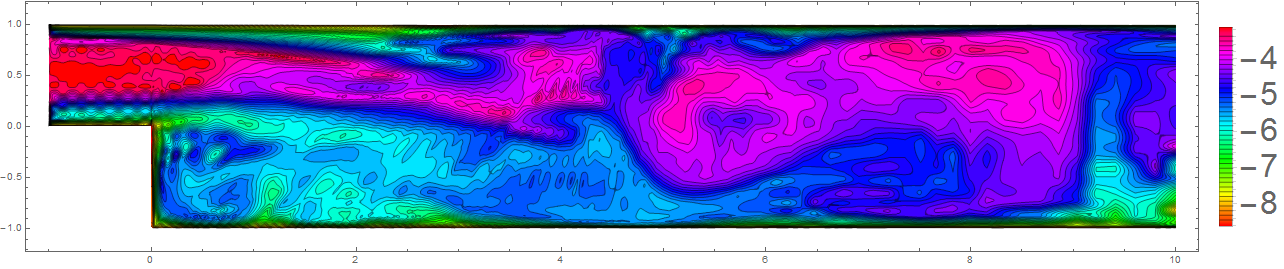}
\end{tabular}
\end{center}
\caption{The solution of turbulent flow based on turbulent model (\ref{eq20}), (\ref{eq21}). The upper pictures illustrate streamlines at the times $2.5s, 5s, 7.5s$ and $10s.$ The lower picture illustrates distribution of the turbulent viscosity (in logaritmic scale) at the time $10s.$}\label{obr:LRN2006}
\end{figure}

\begin{figure}[ht!]
\begin{center}
\includegraphics[width=10.7cm,height=2.7cm]{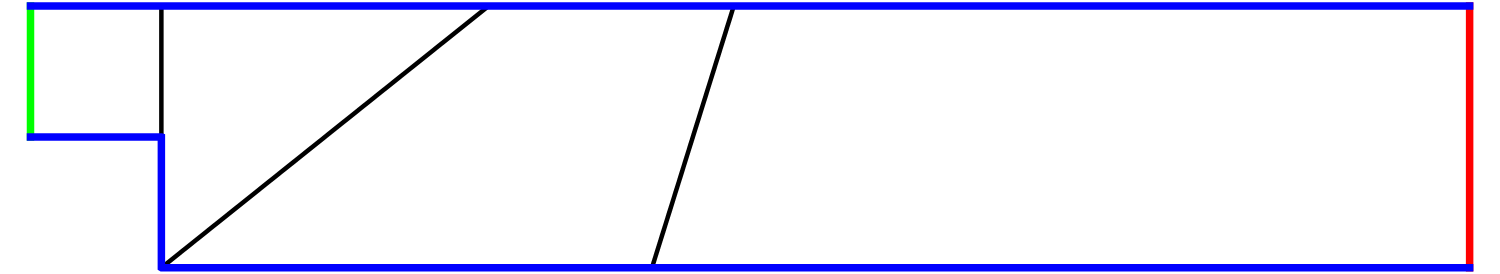}
\end{center}
\caption{L-shaped domain with the NURBS elements.}\label{obr:sit}
\end{figure}

We present the numerical experiment devoted to the turbulent flow with
$Re = 4000$ ($Re = \frac{UL}{\nu},$ $U$  is the maximum velocity of the fluid, $L$ is a characteristic linear dimension, $\nu = 10^{-3}$) through the L-shape domain. All quantities are mentioned in base SI units. The domain with the initial NURBS elements (they are further refined) is shown at the Figure \ref{obr:sit}. We consider the left inflow boundary, right outflow boundary and the remaining boundaries as the solid walls. The velocity $\bm u$ is defined by the parabolic profile with maximum value $\bm u_{\max} = [3, 0]$ at the inflow boundary and zero at the solid walls. The homogeneous Neumann condition is defined at the outflow boundary. The homogeneous Neumann condition is defined for the pressure at the whole boundary. We set $k = 10^{-6}$ and $\omega = 1$ at the inflow boundary and solid walls. At the outflow boundary there is defined the homogeneous Neumann condition too.

The initial condition for the pressure, $k$ and $\omega$ is constant in the whole domain, $p = 0, k = 10^{-6}$ and $\omega = 1.$ The initial velocity distribution is defined by the auxiliary solution of the Stokes problem. The solution of the problems based on turbulent models described in sections \ref{sec31} and \ref{sec32} is illustrated at the Figures \ref{obr:HRN}, \ref{obr:LRN1993} and \ref{obr:LRN2006}. We use time step $\Delta t = 0.05.$

It can be seen the difference between the simulation results obtained by the basic $k-\omega$ model (section \ref{sec31}) and the LRN models (section \ref{sec32}), see Figures \ref{obr:HRN}--\ref{obr:LRN2006}. The main reason is the amount of turbulent viscosity produced by different models in this case, i.e., basic $k-\omega$ model produces more turbulent viscosity than other two models (compare lower pictures in Figures \ref{obr:HRN}--\ref{obr:LRN2006}).

\section{Conclusion}

In the paper, we focused on a numerical modelling of turbulent flows. We
developed and tested an isogeometric analysis based solver for RANS equations with several variants of $k-\omega$ models. The presented results show that the isogeometric analysis is a suitable tool for solving such complex problems. In the future, we plan to generalize our solver to three space dimensions in order to be able to simulate flows in water turbines. Moreover, we will test another turbulent models and suitable stabilization techniques necessary for solving turbulent flow problems with high Reynolds numbers.

\appendix
\section*{Appendix}
The equation (\ref{eq48}) is multiplied by the test function $\varphi_k + \tau_k\bar{\bm u}^u\cdot\nabla\varphi_k$ and the equation (\ref{eq49}) is multiplied by the test function $\varphi_{\omega} + \tau_{\omega}\bar{\bm u}^n\cdot\nabla\varphi_{\omega}.$ Then we integrate over $\Omega$ and use Green's theorem by assuming that $\varphi_k, \varphi_{\omega} \in V_0 = \{\psi \in H^1(\Omega): \psi|_{\partial\Omega} = \bm 0\}.$ We search for $k, \omega \in V = \{\varphi \in H^1(\Omega):\varphi|_{\partial\Omega} = g\}$ so that all functions $\varphi_k \in V_0$ and $\varphi_{\omega} \in V_0$ satisfy
\begin{eqnarray*}
\int_{\Omega}k^{n+1}\varphi_k\mathrm{d}\Omega + \Delta t\int_{\Omega}\left(\frac{\nu_T^n}{\sigma_k} + \nu\right)\nabla k^{n+1}\cdot\nabla\varphi_k\mathrm{d}\Omega + \\
+ \Delta tC_{\mu}\int_{\Omega}\frac{(k^{n+1})^2}{\nu_T^n}\varphi_k\mathrm{d}\Omega + \int_{\Omega}k^{n+1}\tau_k\bar{\bm u}^n\cdot\nabla\varphi_k\mathrm{d}\Omega - \\
- \Delta t\int_{\Omega}\left[\left(\frac{\nu_T^n}{\sigma_k} + \nu\right)\nabla k^{n+1}\right]\cdot\nabla(\tau_k\bar{\bm u}^n\cdot\nabla\varphi_k)\mathrm{d}\Omega + \\
+ \Delta tC_{\mu}\int_{\Omega}\frac{(k^{n+1})^2}{\nu_T^n}\tau_k\bar{\bm u}^n\cdot\nabla\varphi_k\mathrm{d}\Omega = \\
= \int_{\Omega}k^n\varphi_k\mathrm{d}\Omega - \Delta t\int_{\Omega}\bar{\bm u}^n\cdot\nabla k^{n+1}\varphi_k\mathrm{d}\Omega + \Delta t\int_{\Omega}\nu_T^nf^n\varphi_k\mathrm{d}\Omega + \\
+ \int_{\Omega}k^n\tau_k\bar{\bm u}^n\cdot\nabla\varphi_k\mathrm{d}\Omega - \Delta t\int_{\Omega}\bar{\bm u}^n\cdot\nabla k^{n+1}\tau_k\bar{\bm u}^n\cdot\nabla\varphi_k\mathrm{d}\Omega + \\
+ \Delta t\int_{\Omega}\nu_T^nf^n\tau_k\bar{\bm u}^n\cdot\nabla\varphi_k\mathrm{d}\Omega.
\end{eqnarray*}
\begin{eqnarray*}
\int_{\Omega}\omega^{n+1}\varphi_\omega\mathrm{d}\Omega + \Delta t\int_{\Omega}(\sigma_\omega\nu_T^n + \nu)\nabla \omega^{n+1}\cdot\nabla\varphi_\omega\mathrm{d}\Omega + \\
+ \Delta tC_{\omega 2}\int_{\Omega}(\omega^{n+1})^2\varphi_\omega\mathrm{d}\Omega + \int_{\Omega}\omega^{n+1}\tau_\omega\bar{\bm u}^n\cdot\nabla\varphi_\omega\mathrm{d}\Omega - \\
- \Delta t\int_{\Omega}[(\sigma_\omega\nu_T^n + \nu)\nabla \omega^{n+1}]\cdot\nabla(\tau_\omega\bar{\bm u}^n\cdot\nabla\varphi_\omega)\mathrm{d}\Omega + \\
+ \Delta tC_{\omega 2}\int_{\Omega}(\omega^{n+1})^2\tau_\omega\bar{\bm u}^n\cdot\nabla\varphi_\omega\mathrm{d}\Omega = \\
= \int_{\Omega}\omega^n\varphi_\omega\mathrm{d}\Omega - \Delta t\int_{\Omega}\bar{\bm u}^n\cdot\nabla \omega^{n+1}\varphi_\omega\mathrm{d}\Omega + \Delta tC_{\omega 1}\int_{\Omega}f^n\varphi_\omega\mathrm{d}\Omega + \\
+ \int_{\Omega}\omega^n\tau_\omega\bar{\bm u}^n\cdot\nabla\varphi_\omega\mathrm{d}\Omega - \Delta t\int_{\Omega}\bar{\bm u}^n\cdot\nabla \omega^{n+1}\tau_\omega\bar{\bm u}^n\cdot\nabla\varphi_\omega\mathrm{d}\Omega + \\
+ \Delta tC_{\omega 1}\int_{\Omega}f^n\tau_\omega\bar{\bm u}^n\cdot\nabla\varphi_\omega\mathrm{d}\Omega.
\end{eqnarray*}
Now we consider finite dimensional spaces $V^h \subset V$ and $V_0^h \subset V_0$ with basis functions $R_i^k \in V^h$ and $R_i^\omega \in V_0^h.$ The solutions $k_h \in V^h$ and $\omega_h \in V^h$ have the forms of linear combinations of basis functions
\begin{equation}
k_h^{n+1} = \sum\limits_{i=1}^Nk_i^{n+1}R_i^k, \quad k_h^n = \sum\limits_{i=1}^Nk_i^nR_i^k,
\end{equation}
\begin{equation}
\omega_h^{n+1} = \sum\limits_{i=1}^N\omega_i^{n+1}R_i^k, \quad \omega_h^n = \sum\limits_{i=1}^N\omega_i^nR_i^k,
\end{equation}
We consider $\varphi_k \in V_0^h$ and $\varphi_\omega \in V_0^h$. Substituting to the previous equations we have
\begin{eqnarray*}
\sum\limits_{i=1}^Nk_i^{n+1}[\underbrace{\int_\Omega R_i^kR_j^k\textrm{d}\Omega}_{A^k} + \Delta t\underbrace{\int_\Omega\left(\frac{\nu_T^n}{\sigma_k} + \nu\right)\nabla R_i^k\cdot\nabla R_j^k\textrm{d}\Omega}_{C^k} + \\
+ \Delta tC_\mu\underbrace{\int_\Omega\frac{1}{\nu_T^n}R_i^kR_j^k\sum\limits_{l = 1}^Nk_l^{n+1}R_l^k\textrm{d}\Omega}_{B_l^k(k)} + \underbrace{\int_\Omega\left(\tau_k\sum\limits_{m=1}^N\bar{\bm u}_m^nR_m^u\right)\cdot(R_i^k\nabla R_j^k)\textrm{d}\Omega}_{\nabla A^k} +\\
+ \Delta t\underbrace{\int_\Omega(\frac{\nu_T^n}{\sigma_k} + \nu)\nabla R_i^k\tau_k((\sum\limits_{l=1}^N\bar{\bm u}_l^n\cdot\nabla R_l^u)\nabla R_j^k + (\sum\limits_{l=1}^N\bar{\bm u}_l^nR_l^u)(\nabla\cdot\nabla R_j^k))\textrm{d}\Omega}_{\nabla C^k} + \\
+ \Delta t C_\mu\underbrace{\int_\Omega\left(\tau_k\sum\limits_{m=1}^N\bar{\bm u}_m^nR_m^u\right)\frac{1}{\nu_T^n}R_i^k\cdot\left[\nabla R_j^k\left(\sum\limits_{l=1}^Nk_l^{n+1}R_l^k\right)\right]\textrm{d}\Omega}_{\nabla B_l^k(k)} + \\
+ \Delta t\underbrace{\int_\Omega(\sum\limits_{l=1}^N\bar{\bm u}_l^nR_l^u)\cdot(\nabla R_i^kR_j^k)\textrm{d}\Omega}_{D^k} + \\
+ \Delta t\underbrace{\int_\Omega\left(\tau_k\sum\limits_{m=1}^N\bar{\bm u}_m^nR_m^u\right)\cdot\left(\sum\limits_{l=1}^N\bar{\bm u}_l^nR_l^u\right)(\nabla R_i^k\cdot\nabla R_j^k)\textrm{d}\Omega}_{\nabla D^k}] = \\
= \sum\limits_{i=1}^Nk_i^n[\underbrace{\int_\Omega R_i^kR_j^k\textrm{d}\Omega}_{A^k} + \underbrace{\int_\Omega\left(\tau_k\sum\limits_{m=1}^N\bar{\bm u}_m^nR_m^u\right)\cdot(\nabla R_i^kR_j^k)\textrm{d}\Omega}_{\nabla A^k}] + \\
+ \Delta t\underbrace{\int_\Omega\nu_T^nf^nR_j^k\textrm{d}\Omega}_{F^k} + \Delta t\underbrace{\int_\Omega\tau_k(\sum\limits_{l=1}^N\bar{\bm u}_l^nR_l^u)\nu_T^nf^n\cdot\nabla R_j^k\textrm{d}\Omega}_{\nabla F^k},
\end{eqnarray*}
\begin{eqnarray*}
\sum\limits_{i=1}^N\omega_i^{n+1}[\underbrace{\int_\Omega R_i^\omega R_j^\omega\textrm{d}\Omega}_{A^\omega} + \Delta t\underbrace{\int_\Omega(\sigma_\omega\nu_T^n + \nu)\nabla R_i^\omega\cdot\nabla R_j^\omega\textrm{d}\Omega}_{C^\omega} + \\
+ \Delta tC_{\omega 2}\underbrace{\int_\Omega R_i^\omega R_j^\omega(\sum\limits_{l = 1}^N\omega_l^{n+1}R_l^\omega)\textrm{d}\Omega}_{B_l^\omega(\omega)} + \underbrace{\int_\Omega(\tau_\omega\sum\limits_{m=1}^N\bar{\bm u}_m^nR_m^u)(R_i^\omega\nabla R_j^\omega)\textrm{d}\Omega}_{\nabla A^\omega} +\\
+ \Delta t\underbrace{\int_\Omega(\sigma_\omega\nu_T^n + \nu)\nabla R_i^\omega\tau_\omega((\sum\limits_{l=1}^N\bar{\bm u}_l^n\cdot\nabla R_l^u)\nabla R_j^\omega\!\! +\!\! (\sum\limits_{l=1}^N\bar{\bm u}_l^nR_l^u)(\nabla\!\!\cdot\!\!\nabla R_j^\omega))\textrm{d}\Omega}_{\nabla C^\omega} + \\
+ \Delta t C_{\omega 2}\underbrace{\int_\Omega\left(\tau_\omega\sum\limits_{m=1}^N\bar{\bm u}_m^nR_m^u\right)R_i^\omega\cdot\left[\nabla R_j^\omega\left(\sum\limits_{l=1}^N\omega_l^{n+1}R_l^\omega\right)\right]\textrm{d}\Omega}_{\nabla B_l^\omega(\omega)} + \\
+ \Delta t\underbrace{\int_\Omega(\sum\limits_{l=1}^N\bar{\bm u}_l^nR_l^u)\cdot\nabla R_i^\omega R_j^\omega\textrm{d}\Omega}_{D^\omega} + \\
+ \Delta t\underbrace{\int_\Omega\left(\tau_\omega\sum\limits_{m=1}^N\bar{\bm u}_m^nR_m^u\right)\cdot\left(\sum\limits_{l=1}^N\bar{\bm u}_l^nR_l^u\right)(\nabla R_i^\omega\cdot\nabla R_j^\omega)\textrm{d}\Omega}_{\nabla D^\omega}] = \\
= \sum\limits_{i=1}^N\omega_i^n[\underbrace{\int_\Omega R_i^\omega R_j^\omega\textrm{d}\Omega}_{A^\omega} + \underbrace{\int_\Omega\left(\tau_\omega\sum\limits_{m=1}^N\bar{\bm u}_m^nR_m^u\right)R_i^\omega\cdot\nabla R_j^\omega\textrm{d}\Omega}_{\nabla A^\omega}] + \\\
+ \Delta tC_{\omega 1}\underbrace{\int_\Omega f^nR_j^\omega\textrm{d}\Omega}_{F^\omega} + \Delta tC_{\omega 1}\underbrace{\int_\Omega\tau_\omega(\sum\limits_{l=1}^N\bar{\bm u}_l^nR_l^u)f^n\cdot\nabla R_j^\omega\textrm{d}\Omega}_{\nabla F^\omega}.
\end{eqnarray*}

\end{document}